\documentclass[a4paper,14pt]{article}

\usepackage{amssymb}
\usepackage{amsmath}
\usepackage{graphicx}

\newtheorem{thm}{Theorem}[section]

\newtheorem{lem}[thm]{Lemma}
\newtheorem{prop}[thm]{Proposition}
\newtheorem{defn}[thm]{Definition}
\newtheorem{rem}[thm]{Remark}
\numberwithin{equation}{section}

\title{Compactification of the moduli space of $\rho$ vortices}
\author{P. Angulo}

\begin{document}

\maketitle

\begin{abstract}
In this paper we prove the Uhlenbeck compactness of the moduli space of
$\rho$-vortices
introduced in \cite{Ban}, for a k\"{a}hler surface.
\end{abstract}
\section {Introduction}

In this paper, we study the behavior of the set of
solutions to the $\rho$-vortex equations. These equations were introduced
in \cite{Ban}, where it is shown that they feature as specific examples many other
previously studied equations, such as the Hermite-Einstein equations for a
bundle, the vortex equation for a pair and the Seiberg-Witten equations. The
$\rho$-vortex equations correspond to the absolute minima of the Yang-Mills-Higgs
functional (see section \ref{minima of the functional}) when the base space is a
K\"{a}hler manifold. In \cite{Ban}, an analog of the Hitchin-Kobayashi
correspondence is stablished, showing that $\rho$-vortices can be identified
with stable holomorphic pairs, and thus that over a K\"{a}hler manifold, the absolute minima
of the Yang-Mills-Higgs functional can be studied algebraically.

We will prove that the moduli space of $\rho$-vortices
can be given a natural Uhlenbeck compactification, when the base space is a
K\"{a}hler surface. Our argument follows the line of \cite[chapter 4]{DK}, a
standard reference for the proof of the Uhlenbeck compactness of the moduli
space of instantons over a riemannian manifold. All the necessary background
for this paper is covered in \cite{Ban} and this reference (see chapter 2 for
results about connections and a proof of the removal of singularities of Uhlenbeck,
and the appendix for basic results about elliptic operators).

In the particular case $G=GL\left(n,\mathbb{C}\right)$ and $\rho$ the
fundamental representation, the equations are simply called vortex equations.
In this case the construction of the compactification
can be related to that of instantons through dimensional reduction. Using this
relation, the compactification of the space of holomorphic pairs is done in general
dimension in \cite{TiYa}. However, dimensional reduction doesn't seem to work
for a general representation, and a direct approach is necessary.

\section {The $\rho$-vortex equations}

Let $K$ be a compact real Lie group with Lie algebra $\mathfrak{k}$, $G$ the
complexification of $K$, $X$ a compact
k\"{a}hler surface, $P\longrightarrow X$ a principal holomorphic $G$-bundle.
We also have a representation $\rho :K\longrightarrow U\left( V\right) $,
for $V$ a complex vector space which lifts to a representation of $G$ in $%
GL\left( V\right) $. 
We take $G$ reductive, with non-trivial center acting over $V$ as multiplication by scalars, 
and let $E=\rho \left(P\right) $ be the vector bundle $E\times_{\rho}V$.
In section \ref{consequences} we will make some slight hypothesis
over the action. We define an invariant bilinear product in $\mathfrak{g}$ as the killing
form of the semisimple part of the Lie algebra extended to the centre $\mathfrak{z}$. This invariant
metric provides a metric in $Ad(P)$. We also choose a smooth reduction of the
structure group of $P$ to the compact subgroup $K$. When $G=GL\left(
n,\mathbb{C}\right) $, $K=U\left( n\right) $, and $\rho $ is the fundamental
representation, this is equivalent to giving $E$ a bundle metric.

We consider a set of equations for a pair $\left( A,\phi
\right) $ consisting of a connection on the principal bundle $P$ compatible
with the reduction to $K$ (a unitary connection in the case $K=U\left(
n\right) $) and a section of the bundle $E=\rho \left( P\right) $ associated
to $P$ by $\rho $. We will call them the $\rho $-vortex equations.

\begin{equation}\label{rhovortex1}
\Lambda F_{A}-i\rho ^{\ast }(\phi \otimes \phi ^{\ast })+i\tau =0
\end{equation}
\begin{equation}\label{rhovortex2}
F_{A}^{0,2}=F_{A}^{2,0}=0
\end{equation}
\begin{equation}\label{rhovortex3}
\bar{\partial}_{A}\phi=0
\end{equation}

Here $\Lambda $ denotes contraction with the k\"{a}hler form,
$\rho ^{\ast }(\phi \otimes \phi ^{\ast })$ is the unique element of $\mathfrak{g}$
satisfying $\left\langle \rho ^{\ast }(\phi \otimes \phi ^{\ast
}),X\right\rangle =\left\langle \phi ,\rho _{\ast }\left( X\right) \phi
\right\rangle $, for all $X\in g$, $\tau$ is an element in the centre of the Lie algebra
of $K$, $F_{A}^{0,2}$ is the $\left( 0,2\right) $
-part of the curvature in the splitting of 2-forms $\Omega ^{2}=\Omega
^{2,0}\oplus \Omega ^{1,1}\oplus \Omega ^{0,2}$, induced by the complex
structure, and $\bar{\partial}_{A}\phi $ is the $\left( 0,1\right) $-part of
$D_{A}\phi $.

The equations are invariant under the natural action of the gauge group on
$A$ and $\phi$.
$$s\in Gauge_{K}\left( P\right) \qquad s\left( A,\phi \right) =\left(
s^{-1}As,s^{-1}\phi \right) $$
It is natural to consider only solutions modulo gauge equivalence.

A connection is locally expressed as a $\mathfrak{k}$-valued 1-form. The expression of
$A$ in different trivializations differ by the action of the gauge group. In
later sections we will be working with local expressions for the connection
$A$. The following gauge fixing theorem of Uhlenbeck \cite{Uhl:gauge fixing}
will be useful:

\textbf{(Uhlenbeck theorem)}
There is a fixed $\varepsilon $ such that any $g$-valued 1-form $A$ over a
ball $B^{4}$ with $\left\Vert F_{A}\right\Vert \leq \varepsilon $ is gauge
equivalent to $\tilde{A}$, such that $d^{\ast }\tilde{A}=0$. $\tilde{A}$
is said to be in Coulomb gauge. Moreover, $\tilde{A}$
verifies:
$$\|\tilde{A}\|\leq C\|F_{A}\|$$

In a complex vector bundle, given a holomorphic structure and a metric on
the bundle, there is a unique connection, called the Chern connection, which
is compatible with both the holomorphic structure and the metric. In our
principal bundle framework, if we fix a holomorphic structure on the bundle
$P$, then for any reduction of the structure group, there is a unique Chern
connection compatible with both the holomorphic structure and the
reduction. This remark suggest another way of looking at the equations. We
can fix a holomorphic structure on the bundle $P$ and a holomorphic section
$\phi$ of the associated bundle $E$ and ask whether there is a reduction
such that the associated connection solves the $\rho$-vortex equations. The
Hitchin-Kobayashi correspondence asserts that this is the case if and only
if the holomorphic pair $\left( E,\phi \right)$ is $\tau$-stable. Thus the
moduli space of solutions to the equations can be identified with the moduli
space of stable pairs modulo the action of the complex gauge group, which
is the group of automorphisms of the holomorphic pair.

\section{The $YMH$ functional} \label{minima of the functional}

\begin{thm}
The $\rho $-vortices are the absolute minima of the $YMH$ functional:

\begin{equation}\label{functional}
YMH(A,\phi )=\left\| F_{A}\right\| ^{2}
        +\left\| \rho ^{\ast }\left( \phi
            \otimes \phi ^{\ast }\right) -\tau \right\| ^{2}
        +2\left\| d_{\rho A}\left(\phi \right) \right\| ^{2}
\end{equation}
\end{thm}

The solutions of the $\rho $-vortex equations are obviously the
absolute minima of the following functional:
$$YMH(A,\phi )=\int\limits_{X}\left( \left| \Lambda F_{A}-i\rho ^{\ast
}(\phi \otimes \phi ^{\ast })+i\tau \right| ^{2}+4\left| F_{A}^{2,0}\right|
^{2}+4\left| \bar{\partial}_{\rho A}\phi \right| ^{2}\right) $$

We will rewrite this energy functional in a more convenient form.
$$\begin{array}{rl}
\int_{X}\left| \Lambda F_{A}-i\rho ^{\ast }\left( \phi \otimes \phi ^{\ast
}\right) +i\tau \right| ^{2}=&\left\| \Lambda F_{A}\right\| ^{2}+\left\| \rho
^{\ast }\left( \phi \otimes \phi ^{\ast }\right) -\tau \right\| ^{2}+\\
&+2i\int\limits_{X}\left\langle \Lambda F_{A},\tau \right\rangle
-2i\int\limits_{X}\left\langle \Lambda F_{A},\rho ^{\ast }\left( \phi
\otimes \phi ^{\ast }\right) \right\rangle
\end{array}$$

The last two terms admit a more convenient rewriting:
$$\int\limits_{X}\left\langle i\Lambda F_{A},\tau \right\rangle \omega^{2}
  =\int\limits_{X}\left\langle iF_{A},\tau \right\rangle \omega
  =2\pi\deg_{\tau }\left( P\right) $$

The polinomial $\left\langle X,\tau \right\rangle $ is invariant
under conjugation of $X$, so it follows by Chern-Weyl theory that $\deg
_{\tau }\left( P\right) $ is a topological invariant. For the second we
use the k\"{a}hler identities.
$$\begin{array}{rl}
-2i\int\limits_{X}\left\langle \Lambda F_{A},\rho ^{\ast }\left( \phi
\otimes \phi ^{\ast }\right) \right\rangle
&=-2i\int\limits_{X}\left\langle
\Lambda F_{\rho A}\left( \phi \right) ,\phi \right\rangle\\
&=-2i\int\limits_{X}\left\langle \Lambda d_{\rho A}d_{\rho A}\left( \phi
\right) ,\phi \right\rangle\\
&=-2i\int\limits_{X}\left\langle \Lambda \left( \bar{\partial}_{\rho
A}+\partial _{\rho A}\right) \left( \bar{\partial}_{\rho A}+\partial _{\rho
A}\right) \left( \phi \right) ,\phi \right\rangle\\
&=-2i\left( \int\limits_{X}\left\langle i\bar{\partial}_{\rho A}^{\ast }
\bar{\partial}_{\rho A}\left( \phi \right) ,\phi \right\rangle
-\int\limits_{X}\left\langle i\partial _{\rho A}^{\ast }\partial _{\rho
A}\left( \phi \right) ,\phi \right\rangle \right)\\
&=2\left\Vert \partial_{\rho A}\left( \phi\right) \right\Vert
^{2}-2\left\Vert \bar{\partial}_{\rho A}\left( \phi\right) \right\Vert ^{2}
\end{array}$$

We can group the term $\left\| \Lambda F_{A}\right\| ^{2}$ with the other
term in the curvature: $4\left\| F_{A}^{2,0}\right\| ^{2}$ to get $\left\|
\Lambda F_{A}\right\| ^{2}+4\left\| F_{A}^{2,0}\right\| ^{2}=\left\|
F_{A}\right\| ^{2}+8\pi ^{2}Ch_{2}\left( P\right) $, and it follows that the $%
\rho $-vortices are the absolute minima of
$$YMH(A,\phi )=4\pi\deg _{\tau }\left( P\right) +8\pi ^{2}Ch_{2}\left( P\right)
+\left\| F_{A}\right\| ^{2}+\left\| \rho ^{\ast }\left( \phi \otimes \phi
^{\ast }\right) -\tau \right\| ^{2}+2\left\| d_{\rho A}\left( \phi \right)
\right\| ^{2}$$

From now on, we will call the integrand $YMH\left( A,\phi \right)
=\left\vert F_{A}\right\vert ^{2}+\left\vert \rho ^{\ast }\left( \phi
\otimes \phi ^{\ast }\right) -\tau \right\vert ^{2}+2\left\vert d_{\rho
A}\left( \phi \right) \right\vert ^{2}$ the Yang-Mills-Higgs action.

\section{A few consequences of the $\rho $-vortex equations}
\label{consequences}

In this section we will prove some elementary facts about
vortices which will be useful later. We assume that the center $Z$ of the
group $G$ acts on $V$ by multiplication by scalars, so that for any $\xi \in
Z$, $\rho \left( \xi \right) =\lambda(\xi) I$, where $\lambda(\xi)$ is a complex number. 
If we take the inner product of equation \ref{rhovortex1} with $\xi $, multiply by
$\omega ^{2}$, and integrate over $X$, we get:
$$0=\left\langle i\Lambda F_{A},\xi\right\rangle +
\int\left\langle\rho^{\ast}(\phi\otimes\phi^{\ast}),\xi\right\rangle -
\int_{X}\left\langle\tau ,\xi\right\rangle dvol$$
The first and the third terms are constant, so we get:
\begin{equation}\label{norm of phi is constant}
-2\pi deg_{\xi}P+\left\langle \tau,\xi\right\rangle
        vol\left( X\right) =
\int\left\langle\rho^{\ast}(\phi\otimes\phi^{\ast}),\xi\right\rangle
=\int\left\langle\phi\otimes\phi^{\ast},\rho\left( \xi\right) \right\rangle
=\lambda\left(\xi\right) \left\Vert \phi\right\Vert ^{2}
\end{equation}

If we take any $\xi $ such that $\lambda \left( \xi \right) \neq 0$, we
arrive at the conclusion that the $L^{2}$ norm of $\phi $ is constant for
all solutions. If we take $\xi $ in the kernel of $\lambda $ we find that in
order to have solutions, $\tau $ has to verify some linear restrictions. As
$\lambda $ is a linear form, the number of independent restrictions is always
$\dim \left( Z\right) -1$, so the parameter $\tau $ has one degree of
freedom. The formula above implies that the section is zero when $\tau $
satisfies this inequality (for any $\xi $ with $\lambda \left( \xi \right)
>0 $):

$$-2\pi deg_{\xi}P+\left\langle \tau,\xi\right\rangle vol\left(
X\right) \leq0$$

This means that $\tau$ can be put to depend on a real parameter lying in an
interval bounded below. For a more detailed account of the dependence of the
equations on the value of $\tau$ (in the most studied particular cases),
see \cite{decorated bundles}.

A different computation yields a pointwise bound for the norm of the section.
Let $x_{0}$ be the point where $|\phi| $ is maximum. At this point, $\triangle\left|\phi\right|^{2}\leq 0$.
Making use of the compatibility of $\rho\left(A\right)$ with the metric on $E$, we find:

\begin{equation}
\begin{array}{rl}
0\geq \triangle\left|\phi\right|^{2}&=
  \sum_{i}\left(\frac{\partial}{\partial x_{k}}\right)^{2}\left\langle\phi,\phi\right\rangle\\&
  =\sum_{i}\left|\nabla_{k}\phi\right|^{2}+
    \sum_{i}\left\langle \nabla_{k}\nabla_{k}\phi,\phi\right\rangle\\&
  =\left|\nabla\phi\right|^{2}-
    \left\langle D_{A}^{\ast}D_{A}\phi,\phi\right\rangle
\end{array}
\end{equation}

Now we can use the k\"{a}hler identities and the $\rho$-vortex equations to find:
\begin{equation}
0 \geq -\left\langle \partial_{A}^{\ast} \partial_{A} \phi,\phi\right\rangle=
       -\left\langle i\Lambda\bar{\partial}_{A}\partial_{A}\phi,\phi\right\rangle=
       -\left\langle i\Lambda F_{A}\phi,\phi\right\rangle=
       \left\langle \left(\rho ^{\ast }(\phi \otimes \phi ^{\ast })-\tau\right)\phi,\phi\right\rangle
\end{equation}

We can also write it this way:
\begin{equation}
0 \geq \left| \rho ^{\ast }(\phi \otimes \phi ^{\ast })-\tau\right|^{2}+
            \left\langle \rho ^{\ast }(\phi \otimes \phi ^{\ast })-\tau,\tau\right\rangle=
       \left| \rho ^{\ast }(\phi \otimes \phi ^{\ast })-\tau\right|^{2}+  \lambda(\tau)|\phi|^{2}-|\tau|^{2}
\end{equation}

Now if $\lambda\left(\tau\right)> 0$ we have find that
$\left|\phi\right|^{2}\leq\frac{\left|\tau\right|^{2}}{\lambda\left(\tau\right)}$
for the point $x_{0}$ and thus for all points.
If $\lambda\left(\tau\right)\leq 0$, we can use (\ref{bound}) at the end of this section to find:
\begin{equation}\label{}
\left|\phi\right|^{4}\leq
  C\cdot\left| \rho ^{\ast }\left(\phi\otimes \phi ^{\ast }\right) -\tau \right| ^{2}\leq
  C\left(\left|\tau\right|^{2}
    -\lambda\left(\tau\right)\left|\phi\right|^{2}\right)
\end{equation}

Which gives the following bound
\begin{equation}\label{boundphi}
\left|\phi\right|^{2}\leq
    C\max\left\{\left|\tau\right|,1-\lambda\left(\tau\right)\right\}
\end{equation}
This computation can also be made for a group with trivial center, 
showing that in this case there are no solutions to the vortex equations.

The results we prove in this paper fail for $G$ abelian and for actions
which are "partially" abelian.
We will assume from now on that the restriction of the infinitesimal action to the semisimple part of
$\mathbf{g}$ doesn't act trivially over any subspace.

\begin{prop}\label{bound}
\textit{
Let $\rho$ be an action of a reductive group $G$ on a vector
space $V$ satisfying the above hypothesis. For any $\tau \in Re\left(Z\right) $,
there is a constant $C$ such that, for any vector $\phi \in V$:
$$\left\vert \phi\right\vert ^{2}\leq C\left\vert \rho^{\ast}\left(
\phi\otimes\phi^{\ast}\right) -\tau\right\vert $$}
\end{prop}

\textbf{Proof}
We will first show that $\rho ^{\ast }\left( \phi \otimes \phi ^{\ast
}\right) +\tau =0$ implies $\phi =0$. We take an element $X$ of the
semisimple part of $\mathbf{k}$:
$$0=\left\langle \rho^{\ast}\left( \phi\otimes\phi^{\ast}\right)
-\tau,X\right\rangle =\left\langle \phi,\rho_{\ast}\left( X\right)
\phi\right\rangle $$

Exponentiating we find:
$$\left\langle \phi ,\rho _{\ast }\left( g\right) \phi \right\rangle =
    \left|\phi \right| ^{2}
\Rightarrow \rho _{\ast }\left( g\right) \phi
    =e^{i\theta }\phi\quad \forall g \in K$$

This means that the semisimple part of $K$ (and thus also
that of $G$) acts in an abelian way. We deduce that
$e^{i\theta }=1$, which by hypothesis implies $\phi =0$.

If the thesis didn't hold, we could find a sequence of $\phi_{i}$ with
$\frac{\left\vert \phi_{i}\right\vert ^{2}}{\left\vert \rho^{\ast}\left(
\phi_{i}\otimes\phi_{i}^{\ast}\right) +\tau\right\vert }$ tending to $\infty$. We
first notice that the norms $\left\vert \phi_{i}\right\vert $ are bounded
below if $\tau\neq0$, because if $\left\vert \phi\right\vert ^{2}\leq
\min\left\{ \frac{\left\vert \tau\right\vert }{2C},\frac{\left\vert
\tau\right\vert }{2}\right\} $ we get
$$\left\vert \rho^{\ast}\left( \phi\otimes\phi^{\ast}\right) +\tau\right\vert
\geq\left\vert \tau\right\vert -\left\vert \rho^{\ast}\left( \phi\otimes
\phi^{\ast}\right) \right\vert \geq\left\vert \tau\right\vert -C\left\vert
\phi\right\vert ^{2}$$

so $\left\vert \phi\right\vert ^{2}\leq\frac{\left\vert \tau\right\vert }{2}
\leq\left\vert \tau\right\vert -C\left\vert \phi\right\vert
^{2}\leq\left\vert \rho^{\ast}\left( \phi\otimes\phi^{\ast}\right)
+\tau\right\vert $ and $\frac{\left\vert \phi\right\vert ^{2}}{\left\vert
\rho^{\ast}\left( \phi\otimes\phi^{\ast}\right) +\tau\right\vert }\leq1$

So if $\tau\neq0$ we can take a subsequence such that the sequence of norms $%
\left\vert \phi_{i}\right\vert $ has a limit which is either $k>0$ or $%
\infty $. Then we consider the vectors $\varphi_{i}=\frac{\phi_{i}}{%
\left\vert \phi_{i}\right\vert }$, and extracting another subsequence, we
can assume that the vectors $\varphi_{i}$ converge to a unit vector $\varphi$
with $\rho ^{\ast}\left( \varphi\otimes\varphi{}^{\ast}\right) +a=0$, where $%
a$ is either $\frac{\tau}{k}$ or zero, where we use the preceding remark for
the case $\tau\neq0$. Both of this possibilities contradict the preceding
paragraph, and so the proof is completed.

\section{Local theory} \label{local}

\begin{defn}
Given an open subset $U$ of the k\"{a}hler manifold $X$ homeomorphic to
$B^{4}$, a pair $\left(A,\phi\right)$ consisting of a $\mathfrak{g}$-valued 1-form $A$
on $X$, together with a $V$-valued function $\phi $ on $X$ is called a local
$\rho $-pair. If we say a pair is $L_{k}^{2}$, we mean $A\in L_{k}^{2}\left(
U,\Omega ^{1}\otimes \mathfrak{g}\right) $, $\phi \in L_{k}^{2}\left( U,V\right) $.
\end{defn}

\begin{defn}
A local $\rho $-pair is a local $\rho $-vortex iff:
\begin{enumerate}
\item $d^{\ast }A=0$
\item $d^{+}A+\left(A\wedge A\right)^{+}
-i\omega \rho ^{\ast }\left( \phi \otimes \phi
^{\ast }\right) +i\omega \tau =0$
\item $\bar{\partial}\phi +\rho \left( A^{0,1}\right) \phi =0$
\end{enumerate}
This equations are called the local $\rho$-vortex equations.
\end{defn}

A local $\rho $-pair is the local version of the $\rho $-vortex equations.
Note that we make use of the decomposition $\Omega ^{+}=\Omega ^{2,0}\oplus
\Omega ^{0,2}\oplus \Omega ^{0}\omega $ to put together all the equations
involving the curvature of the connection $A$. Using the gauge fixing theorem
of Uhlenbeck, we note that any $\rho $-vortex gives rise to a local $\rho$
-vortex over any trivializing open set, just by taking the Coulomb gauge
representative with respect to the trivial connection. The gauge fixing
theorem also asserts that $\left\Vert A\right\Vert _{L_{1}^{2}\left(
U\right) }\leq C\left\Vert F_{A}\right\Vert _{L^{2}\left( U\right) }$, or,
using lemma \ref{bound}, that

\begin{equation}\label{local_bound}
\left\Vert A\right\Vert _{L_{1}^{2}\left( U\right)}+
\left\Vert \phi \right\Vert _{L^{4}\left( U\right) }
\leq C\int_{U} YMH\left(A,\phi \right)
\end{equation}

\begin{lem} \label{keylemma}
There is an $\varepsilon $ such that any local $\rho $-vortex over an
open subset $U$ of $X$ with $\left\Vert A\right\Vert _{L_{1}^{2}\left(
U\right) }+\left\Vert \phi \right\Vert _{L^{4}\left( U\right) }+vol\left(
U\right) \leq \varepsilon $ verifies:
$$\left\Vert A\right\Vert _{L_{p}^{2}\left( V\right) }+\left\Vert
\phi\right\Vert _{L_{p}^{2}\left( V\right) }\leq M_{p}\left( \left\Vert
A\right\Vert _{L_{1}^{2}\left( U\right) }+\left\Vert \phi\right\Vert
_{L^{4}\left( U\right) }+vol\left( U\right) \right) $$

for any interior domain $V\Subset U$ \ and $p\geq 1$. The number $M$ depend
on $p$ and $V$.
\end{lem}

\textbf{Proof}. We will try to use the bootstrapping technique with the $\rho $
-vortex equations to bound all the Sobolev norms. We first cut off $A$ and 
$\phi $ so that they are zero at the boundary. Let $U_{1}$ be such that
$V\Subset U_{1}\Subset U$, and let $\psi $ be a cut off
function which is one over $U_{1}$ and $0$ away from $U$.
We will use the basic facts about elliptic operators which can be found in \cite{We},
together with the Sobolev embedding and multiplication theorems, which we recall in the appendix.

The over-ellipticity of $\bar{\partial}$, combined with the fact that
$\psi\phi$ is zero on the boundary, yields:
$$\begin{array}{rl}
\left\| \psi \phi \right\| _{L_{1}^{2}\left( U\right) }&
\leq C\left( \left\| \bar{\partial}\left( \psi \phi \right) \right\|
_{L^{2}\left( U\right) }+\left\| \psi \phi \right\| _{L^{2}\left( U\right)
}\right)\\&
\leq C\left( \left\| A\right\| _{L_{1}^{2}\left( U\right) }\left\|
\psi \phi \right\| _{L_{1}^{2}\left( U\right) }+\left\| \phi \right\|
_{L^{2}\left( U\right) }\right)
\end{array}
$$

So if $C\varepsilon\leq\frac{1}{2}$ we can rearrange the inequality:
$$\begin{array}{rl}
\left\| \phi \right\| _{L_{1}^{2}\left( U_{1}\right) }&
\leq \left\| \psi\phi \right\| _{L_{1}^{2}\left( U\right) }\\&
\leq M_{0}\left\| \phi \right\|_{L^{2}\left( U\right) }
\end{array}$$

Thus there is a bound for the $L_{1}^{2}$-norm of $\phi $ over the smaller
open set $U_{1}$. Let's say $\left\| A\right\| _{L_{1}^{2}\left(
U_{1}\right) }+\left\| \phi \right\| _{L_{1}^{2}\left( U_{1}\right) }\leq
\delta $. In order to bound the $L_{2}^{2}$ norms we repeat the above
construction, taking another open subset $U_{2}$, and a cut off function $%
\psi $ which is one over $U_{2}$ and zero out of $U_{1}$:
$$\begin{array}{rl}
\left\| \psi A\right\| _{L_{2}^{2}\left( U_{1}\right) }&
\leq C\left( \left\|\left( d^{\ast }+d^{+}\right) \left( \psi A\right) \right\|
_{L_{1}^{2}\left( U_{1}\right) }+\left\| \psi A\right\| _{L^{2}\left(
U_{1}\right) }\right)\\&
\leq C\{\left\| d\psi \right\| _{L_{3}^{2}}\left\| A\right\|
_{L_{1}^{2}\left( U_{1}\right) }+\left\| \psi A\right\| _{L_{2}^{2}\left(
U_{1}\right) }\left\| A\right\| _{L_{1}^{2}\left( U_{1}\right) }+\\&
\qquad+\left\| \psi \phi \right\| _{L_{2}^{2}\left( U_{1}\right) }\left\|
\phi \right\| _{L_{1}^{2}\left( U_{1}\right) }+\tau ^{2}\int \psi +\\&
\qquad+\left\|\psi \right\| _{L_{3}^{2}\left( U_{1}\right) }\left\| A\right\|
_{L^{2}\left( U_{1}\right) }\}\\
\left\Vert \psi\phi\right\Vert _{L_{2}^{2}\left( U_{1}\right) }
&\leq C\left(\left\Vert \bar{\partial}\left( \psi\phi\right) \right\Vert
_{L_{1}^{2}\left( U_{1}\right) }+\left\Vert \psi\phi\right\Vert
_{L^{2}\left( U_{1}\right) }\right)\\&
\leq C( \left\| d\psi \right\| _{L_{3}^{2}}\left\| \phi
\right\| _{L_{1}^{2}\left( U_{1}\right) }+\left\| \psi A\right\|
_{L_{2}^{2}\left( U_{1}\right) }\left\| \phi \right\| _{L_{1}^{2}\left(
U_{1}\right) }+\\&
\qquad+\left\| A\right\| _{L_{1}^{2}\left( U_{1}\right) }\left\|
\psi \phi \right\| _{L_{2}^{2}\left( U_{1}\right) }+\left\| \psi \phi
\right\| _{L^{2}\left( U_{1}\right) })
\end{array}$$

If $C\delta<\frac{1}{2}$\ we can rearrange the second inequality:
$$\left\Vert \psi\phi\right\Vert _{L_{2}^{2}\left( U_{1}\right) }
\leq C\left(\left\Vert \psi A\right\Vert _{L_{2}^{2}\left( U_{1}\right) }
\left\Vert\phi\right\Vert _{L_{1}^{2}\left( U_{1}\right) }+\left\Vert \phi\right\Vert
_{L_{1}^{2}\left( U_{1}\right) }\right) $$

and inserting this into the first one and rearranging (if $\delta$
is small):
$$\begin{array}{rl}
\left\| \psi A\right\| _{L_{2}^{2}\left( U_{1}\right) }&
\leq M\{\left\|A\right\| _{L_{1}^{2}\left( U_{1}\right) }+
\left\| \psi A\right\|_{L_{2}^{2}\left( U_{1}\right) }
\left\| A\right\| _{L_{1}^{2}\left(U_{1}\right) }+\\&
+\left\| \psi A\right\| _{L_{2}^{2}\left( U_{1}\right) }\left\| \phi
\right\| _{L_{1}^{2}\left( U_{1}\right) }^{2}+\left\| \phi \right\|
_{L_{1}^{2}\left( U_{1}\right) }^{2}+\tau ^{2}vol\left( U_{2}\right) \}
\end{array}$$
$$\left\| A\right\| _{L_{2}^{2}\left( U_{2}\right) }\leq \left\| \psi
A\right\| _{L_{2}^{2}\left( U_{1}\right) }\leq M\left( \left\| A\right\|
_{L_{1}^{2}\left( U_{1}\right) }+\left\| \phi \right\| _{L^{2}\left(
U_{1}\right) }+vol\left( U_{1}\right) \right) $$

The same argument provides a bound for $\left\Vert A\right\Vert
_{L_{3}^{2}\left( U_{3}\right) }+\left\Vert \phi\right\Vert
_{L_{3}^{2}\left( U_{3}\right) }$, but in order to rearrange the inequality,
we must demand that $\left\Vert A\right\Vert _{L_{2}^{2}\left( U_{2}\right)
} $ is small, and this amounts to a bound on $\left\Vert A\right\Vert
_{L_{1}^{2}\left( U\right) }+\left\Vert \phi\right\Vert _{L^{4}\left(
U\right) }+vol\left( U\right) $.

For higher $k$ it is no longer necessary to rearrange the inequalities, so
that no further restriction needs to be imposed on $\varepsilon $. We prove
the result inductively for all $k$.
$$\begin{array}{rl}
\left\| A\right\| _{L_{k+1}^{2}\left( U_{k+1}\right) }&
\leq C\left\| \left(d^{+}+d^{\ast }\right) \left( A\right) \right\|
 _{L_{k}^{2}\left(U_{k}\right) }\\&
\leq C\left( \left\| A\right\| _{L_{k}^{2}\left( U_{k}\right)
}^{2}+\left\| \phi \right\| _{L_{k}^{2}\left( U_{k}\right) }^{2}+\tau
^{2}vol\left( U_{k+1}\right) \right)\\&
\leq M\left( \left\|A\right\| _{L_{k}^{2}\left( U_{k}\right) }
+\left\| \phi \right\|_{L_{k}^{2}\left( U_{k}\right) }+
vol\left( U_{k+1}\right) \right)\\&
\leq M\left( \left\| A\right\| _{L_{1}^{2}\left( U\right) }+\left\|
\phi \right\| _{L^{2}\left( U\right) }+vol\left( U^{\prime }\right) \right)
\end{array} $$

Here $M$ is a bound on $\left( \left\| A\right\| _{L_{k}^{2}\left(
U_{k}\right) }+\left\| \phi \right\| _{L_{k}^{2}\left( U_{k}\right)
}+vol\left( U_{k+1}\right) \right) $. A bound for $\left\| \phi \right\| $
can be found analogously. Along the proof we obtain bounds for the norms
over successively smaller open sets, but this can be done so that they all
contain the fixed open set $V$.

\begin{rem} \label{regularity}
The above proof shows that a $L_{1}^{2}$ local pair over $U$ with small
action is actually $C^{\infty }$. \ If we start with a $L_{1}^{2}$ local
pair we find that all its Sobolev norms are bounded so the pair is in $%
L_{k}^{2}$ for all $k$, and the Sobolev embedding theorem shows that the
pair is $C^{\infty }$.
\end{rem}

\section{The limit pairs}

In this section we find the desired compactification of the moduli space of
$\rho $-vortices. Given a sequence of $\rho $-vortices
$\left( A_{k},\phi _{k}\right) $, the action densities
$YMH\left( A_{k},\phi _{k}\right) $ can be though of as bounded measures
over $X$. Then we can take a subsequence such that they converge weakly
to a measure $\nu $ in $\mathcal{M}\left( X\right) $.

Let's now consider a point in $X$ which have a neighborhood $U$ of $\nu $
measure less than the $\varepsilon $ of \ref{keylemma}. We put all the
connections in a local Coulomb gauge. Using \ref{local_bound}, we see that
for sufficiently big $k$:

$$\left\Vert A\right\Vert _{L_{1}^{2}\left( U\right)}+
\left\Vert \phi \right\Vert _{L^{4}\left( U\right) }
\leq\int_{U} YMH\left( A_{k},\phi_{k}\right)
\leq \varepsilon $$

Then we can apply lemma (\ref{keylemma}) and find that all the $L_{p}^{2}\left( U\right) $
Sobolev norms of $A$ and $\phi $ are bounded by numbers independent of $k$,
and extract a subsequence converging in $C^{\infty }$ to a solution $\left(
A,\phi \right) $, so that the measure $\nu $ is actually equal to
$\left| F_{A}\right| ^{2}+$ $\left| \rho ^{\ast }\left( \phi \otimes \phi
^{\ast }\right) -\tau \right| ^{2}+$ $2\left| d_{\rho A}\left( \phi \right)
\right| ^{2}$ in that neighborhood. The diagonal argument of \cite{DK} allows to
find a single gauge and subsequence for which $\left( A_{k},\phi _{k}\right)$
converge in $C^{\infty }\left( U\right) $ to a $C^{\infty }\left( U\right)$
solution in all the points with a good neighborhood.

As $YMH_{X}\left( A_{k},\phi _{k}\right) $ is a constant $C$ (recall section
\ref{functional}), we see that $\nu \left( X\right) =C$ and there can be at most
$\left\lfloor \frac{C}{\varepsilon }\right\rfloor +1$ points $x_{k}$ without
a neighborhood of $\nu $ measure less than $\varepsilon $, for in the
contrary we could take disjoint balls around $\left\lfloor \frac{C}{
\varepsilon }\right\rfloor +1$ points, and the $\nu $ measure of $X$ would
be greater than $C$. Then the measure $\nu $ can be identified with a
continuos measure $\left| F_{A}\right| ^{2}+$ $\left| \rho ^{\ast
}\left( \phi \otimes \phi ^{\ast }\right) -\tau \right| ^{2}+$ $2\left|
d_{\rho A}\left( \phi \right) \right| ^{2}$ outside the points $x_{k}$, and
we have:
\begin{equation} \label{measure}
\nu =\left| F_{A}\right| ^{2}+\left| \rho ^{\ast }\left( \phi
\otimes \phi ^{\ast }\right) -\tau \right| ^{2}+2\left| d_{\rho A}\left(
\phi \right) \right| ^{2}+\sum\limits_{i}8\pi ^{2}n_{i}\delta _{x_{i}}
\end{equation}

The pair $\left( A,\phi \right) $ is only defined in $X\backslash
\{x_{1},..x_{k}\}$. In the next section we will show that the pair can be extended
to another pair defined on the whole manifold $X$. 
This limit pair is defined over a bundle $\tilde{P}$
which can be different to the original $P$. The result is analogous to the
removal of singularities of Uhlenbeck.
In the proof we will use that $\phi$ is uniformly bounded in $X\backslash
\{x_{1},..x_{k}\}$. This follows from \ref{boundphi}.

The next step is to show that the $n_{i}$ are integers.
For this we apply the Chern Weil theory. 
The total action of a $\rho $-vortex is a topological invariant,
namely $4\pi\deg _{\tau }\left( P\right) +8\pi ^{2}Ch_{2}\left( P\right) $. The
degree of two bundles which are isomorphic in $X\backslash \Gamma $ for
$\Gamma $ a finite set of points is the same, but $Ch_{2}\left( P\right) $ may
change. However, $Ch_{2}\left( P\right) $ is integer valued, and thus
$\dfrac{1}{8\pi ^{2}}\int_{B_{r}\left( x_{i}\right) }YMH\left( A,\phi \right)
mod \,\mathbb{Z}$ depends only on the values of the pair on the boundary, because any
two different pairs which are equal in the boundary can be continued in the same
way to the whole $X$ to give total $YMH\left( A,\phi \right)$ actions whose difference is
$8\pi$-times an integer.
In a neighborhood of one of the $x_{i}$, we can consider all the pairs $\left(
A_{k},\phi _{k}\right) $, and the limiting pair $\left( A,\phi \right) $,
extended through the singularity. In a fixed small sphere around $x_{i}$,
$\left( A_{k},\phi _{k}\right) $ are smooth $\rho $-vortices converging in
$C^{\infty }$ to $\left( A,\phi \right) $, and we deduce that
$\lim_{k\rightarrow \infty }\int_{B_{r}\left( x_{i}\right) }YMH\left(
A_{k},\phi _{k}\right) =\int_{B_{r}\left( x_{i}\right) }YMH\left( A,\phi
\right) +8\pi ^{2}n_{i}$, with $n_{i}$ an integer. Taking $r$
small enough, we see that this $n_{i}$ is the one appearing in (\ref{measure}),
and we also see that this integer is positive.

The convergence is in $C^{\infty }$ away from the points $x_{i}$, and this
means that the total action $\int_{X}YMH\left( A,\phi \right) $ is equal to
$\lim\int_{X}YMH\left( A_{k},\phi _{k}\right) -\sum 8\pi ^{2}n_{i}$. This allows
the calculation of $Ch_{2}\left( \tilde{P}\right) $. Now it is possible to
define the Uhlenbeck compactification of the moduli of stable $\rho $-pairs:
the moduli of solutions $M_{d,n}$ with a bundle $P$ with $\deg _{\tau }P=d$
and $Ch_{2}\left( P\right) =n$ is embedded into the space of ideal pairs:
$$\tilde{M}_{d,n}=\amalg _{k}M_{d,n-k}\times Sym^{k}\left( X\right) $$
where a sequence $\left( \left( A,\phi \right)_{k} ,\bar{x}_{k}\right)$
converges to $\left( \left( A,\phi \right) ,\bar{x}\right) $ whenever $\bar{x}
=\left( \bar{y},\bar{z}\right) $, $\bar{x}_{k}$ converges to $\bar{y}$, the
pair $\left( A,\phi \right) _{k}$ converges in $C^{\infty }\left(
X\backslash \{\bar{z}\}\right) $ to $\left( A,\phi \right) $, and the
measures $YMH\left( A,\phi \right) _{k}+\sum 8\pi ^{2}n_{i}\delta
_{x_{k}^{i}}$ converge weakly to $YMH\left( A,\phi \right) +\sum 8\pi
^{2}n_{i}\delta _{x^{i}}$. The closure of $M_{d,n}$ inside $\tilde{M}_{d,n}$
is the desired compactification.

\section{Removal of singularities}
\begin{thm}
\textit{Let $\left( A,\phi \right) $ be a smooth $\rho $-vortex on a bundle over
$B^{4}\backslash \left\{ 0\right\} $ (with respect to a given k\"{a}hler
structure on $B^{4}$) satisfying:}
$$\left\Vert F_{A}\right\Vert ^{2}
+\frac{1}{4}\left\Vert \rho ^{\ast}\left( \phi \otimes \phi ^{\ast }\right)-\tau \right\Vert ^{2}
+\left\Vert D_{\rho A}\left( \phi \right) \right\Vert ^{2}<\infty $$
$$\left\Vert\phi\right\Vert_{L^{\infty}}<\infty$$

\textit{Then there is a $\rho $-vortex on $B^{4}$ whose restriction to
$B^{4}\backslash \left\{ 0\right\} $ is gauge equivalent to $\left( A,\phi
\right) $.}
\end{thm}

\textbf{Proof}.
In order to construct a pair on $B^{4}$ from a pair on
$B^{4}\backslash \left\{ 0\right\} $, we follow quite the same procedure as
in the fourth chapter of \cite{DK}. 
Let $D\left( r\right) $ be the ball of radius $r$ around $0$,
$N\left( r\right) =D\left( 4r\right) \backslash D\left( r\right) $,
$N^{\prime }\left( r\right) =D\left( 3r\right) \backslash D\left( 2r\right) $
, and $\psi _{r}$ a cut off function which is $0$ over $D\left( 2r\right) $
and $1$ outside $D\left( 3r\right) $. Then we construct $\psi _{r}A$ in the
same way as in \cite[p. 168]{DK}, where it is shown that $\left\Vert F_{\psi
_{r}A}^{+}\right\Vert _{L^{2}\left( N^{\prime }\left( r\right) \right) }\leq
C\left\Vert F_{A}\right\Vert _{L^{2}\left( N\left( r\right) \right) }$. We
also let $\psi _{r}\phi $ be the product of the section $\phi $ with the cut
off function, extended by zero. In this way the pair $\left( A,\phi \right) $
is extended to a pair defined over $B^{4}$ which is no longer a $\rho $
-vortex, but can be shown to have small energy:

\begin{equation}\label{bound on the action}
\begin{array}{l}
\qquad\int\limits_{B^{4}}\left( \left\vert
\Lambda F_{\psi _{r}A}-i\rho ^{\ast }(\psi _{r}\phi \otimes \psi _{r}\phi
^{\ast })+\tau \right\vert ^{2}+4\left\vert F_{\psi _{r}A}^{0,2}\right\vert
^{2}+4\left\vert \bar{\partial}_{\rho \psi _{r}A}\left( \psi _{r}\phi
\right) \right\vert ^{2}\right)\\
=\int\limits_{N^{\prime }\left( r\right) }( \left\vert \Lambda
F_{\psi _{r}A}- i\rho ^{\ast }(\psi _{r}\phi \otimes \psi _{r}\phi
^{\ast })+\tau \right\vert ^{2}+4\left\vert F_{\psi
_{r}A}^{0,2}\right\vert ^{2}+4\left\vert \bar{\partial}_{\rho \psi
_{r}A}\left( \psi _{r}\phi \right) \right\vert ^{2}\\
\qquad\qquad\qquad +\tau vol\left(D\left( 2r\right) \right))\\
\leq\int \limits_{N^{\prime}\left( r\right) }\left\vert
F_{\psi_{r}A}^{+}\right\vert ^{2}+\int \limits_{N^{\prime}\left( r\right)}
\left(\left\vert \rho^{\ast}
\left( \psi_{r}\phi\otimes\psi_{r}\phi^{\ast}\right)
\right\vert ^{2}+4\left\vert \bar{
\partial}_{\rho\psi_{r}A}\left( \psi_{r}\phi\right) \right\vert ^{2}\right)+\tau^{2}
vol\left( D\left( 2r\right) \right)\\
\leq C\left( \left\Vert F_{A}\right\Vert _{L^{2}\left( N\left(
r\right) \right) }^{2}+\int \limits_{N^{\prime}\left( r\right) }\left(\left\vert
\rho^{\ast}\left( \phi\otimes\phi^{\ast}\right) \right\vert ^{2}+4\left\vert
\bar{\partial}_{\rho A}\phi\right\vert ^{2}+\left\vert \phi\right\vert ^{2}
\right)+r\right)
\end{array}
\end{equation}
We notice the RHS tends to zero when $r\rightarrow 0$.

If we take $R$ small enough, we can assume that
$\left\Vert F_{A}\right\Vert _{L^{2}(B(R))}$ is as small as we need.
We define for this $R$ a cutoff function $\psi$ which is one inside $B(R)$ and zero
outside $B(2R)$.
Let $r_{i}$ be a sequence of positive radii converging to $0$.
Then $(A\psi_{r}\psi,\phi\psi_{r}\psi)$ defines a sequence of pairs on
the trivial bundle over $S^{4}$ after extending them by zero.
The only point in moving our pairs into $S^{4}$ is that we don't have to care about
boundary conditions when using the Sobolev inequalities.
For each $i$,
we can use the Uhlenbeck gauge fixing theorem to find a local representative $(A_{i},\phi_{i})$
for the pair $(A_{r_{i}},\phi_{r_{i}})$
such that $d^{\ast}A_{i}=0$ and $\|A_{i}\|_{L^{2}_{1}}\leq \varepsilon$. We also have $\|\phi_{i}\|_{L^{2}_{1}}\leq C$
and so we can extract a weak $L^{2}_{1}\oplus L^{2}_{1} \, (B^{4})$
convergent subsequence, and the limit $(A_{\infty},\phi_{\infty})$
is a $\rho$-vortex inside $B(R)$ by (\ref{bound on the action}).

Over any compact subset of $B(R)\backslash\{0\}$, the pairs
$(A_{i},\phi_{i})$ are eventually gauge equivalent to
$(A,\phi)$. It follows as in \cite{DK} that the limit is gauge equivalent to $(A,\phi)$
over all $B(R)\backslash\{0\}$.

It still remains to show that the limit pair is smooth in a neighborhood of $0$.
We will first simply try to do some bootstrapping with the equations. We know that
$A\in L^{2}_{1}\subset L^{4}$ and $\phi\in L^{\infty}$. The equation
$$\bar{\partial}\phi+\rho(A)\phi=0 $$
allows us to bound $\bar{\partial}\phi$ over $B(R)$, but, as $\bar{\partial}$ is an overdetermined
elliptic operator, we see that
$$\|\phi\|_{L^{4}_{1}}\leq C\|\rho(A)\phi\|_{L^{4}}\leq
  C\|A\|_{L^{2}_{1}}\|\phi\|_{L^{\infty}}
  \leq \infty $$
and using the same equation again:
$$\|\phi\|_{L^{2}_{2}}\leq C\|\rho(A)\phi\|_{L^{2}_{1}}\leq
C\left( \|A\|_{L^{2}_{1}}\|\phi\|_{L^{4}_{1}}
 \right)\leq \infty$$
where we have used the Sobolev multiplication theorem on the appendix.
We deduce that $\|\phi\|_{L^{2}_{2}}$ is bounded in $S^{4}$.

This method will not give any further bound for the Sobolev norms of $A$ or $\phi$.
But we can show that $A$ is in $L^{2}_{2}(S^{4})$ using a slightly modified version of \cite[4.4.13]{DK}.
We think of $\phi$ as fixed, and study $A$ as an $L^{2}_{1}$ solution of the following:

\begin{equation}\label{eqsforA}
\begin{array}{c}
\|A\|_{L^{2}_{1}}\leq \varepsilon\\
d^{\ast}A=0\\
d^{+}A+(A\wedge A)^{+}
 -\frac{i}{2}\left(\rho ^{\ast }(\phi \otimes \phi ^{\ast })+\tau\right)\omega
\textrm{ is smooth}.
\end{array}
\end{equation}

We can also state the last condition by saying that $d^{+}A+(A\wedge A)^{+}=\alpha$,
for some $\alpha$ in $L^{2}_{2}$.
This is so because $\phi\in L^{2}_{2}$, and thus
$(\phi \otimes \phi ^{\ast })\in L^{2}_{2}$.
We can assume that $\|\alpha\|_{L^{2}_{2}} $ is small, if $R$ was chosen
adequately.
We can embed equation (\ref{eqsforA}) into a 1-parameter family of equations:
\begin{equation}\label{*t1}
\|B\|_{L^{2}_{1}}\leq \varepsilon
\end{equation}
\begin{equation}\label{*t2}
d^{\ast}B=0
\end{equation}
\begin{equation}\label{*t3}
d^{+}B+(B\wedge B)^{+}=t\alpha
\end{equation}
and use the method of continuity to find a $L^{2}_{2}$ solution to the equations above.
We can easily get a priori bounds for any solution of the equations.
If $B$ is a $L^{2}_{2}$ solution for some $t$, we have:

$$\|B\|_{L^{2}_{2}}
\leq C\|(d^{\ast}+d^{+})B\|_{L^{2}_{1}}
\leq C\left(\|B\|_{L^{2}_{2}}\|B\|_{L^{2}_{1}}+\|\alpha\|_{L^{2}_{1}}\right)$$
and then we can rearrange the inequality to get the following bound:
$$\|B\|_{L^{2}_{2}}\leq C\|\alpha\|_{L^{2}_{1}}$$
And also:
$$\|B\|_{L^{2}_{3}}
\leq C\|(d^{\ast}+d^{+})B\|_{L^{2}_{2}}
\leq C\|B\|_{L^{2}_{2}}^{2}+\|\alpha\|_{L^{2}_{2}}
\leq\infty$$

If we have a sequence of times $t_{i}\rightarrow t$, we can extract a $L^{2}_{2}$ strongly convergent
subsequence, and the limit must be a solution for time $t$.

If we know that there is a solution for a time $t$, we can use the inverse function theorem
to show that there are solutions to (\ref{*t2}) and (\ref{*t3}) for times $t'$ close to $t$, but we
can only achieve $\|B_{t'}\|_{L^{2}_{1}}\leq 2\varepsilon$. However, for any solution, we have:
$$\|B\|_{L^{2}_{1}}
\leq C\|(d^{\ast}+d^{+})B\|_{L^{2}}
\leq C\left(\|B\|_{L^{2}_{1}}^{2}+\|\alpha\|_{L^{4}}^{2}\right)$$
and thus if $\|B\|_{L^{2}_{1}}\leq 2\varepsilon$, and $\|\alpha\|_{L^{4}}$ is small, we arrive at (\ref{*t1}).

Thus the set of times for which there is a solution is both open and close, and non-empty, and so there
is a solution for all times. In particular, we have found a $L^{2}_{2}$ solution for $t=1$.
If there are two $L^{2}_{1}$ solutions $A$ and $B$ to (\ref{eqsforA}),
we can take their difference and compute its norm:

$$\|A-B\|_{L^{2}_{1}}
\leq C\|(d^{\ast}+d^{+})(A-B)\|_{L^{2}}
\leq C\|B+A\|_{L^{2}_{1}}\|A-B\|_{L^{2}_{1}}$$

If $\|B+A\|_{L^{2}_{1}}$ is small, we must have $A=B$.
This completes the proof that $A\in L^{2}_{2}(S^{4})$. The proof that the pair is smooth now follows by the
standard elliptic bootstrapping, combined with the Sobolev multiplication theorem.

\section{Appendix}
Let $L^{p}_{k}$ be the space of measurable functions over a compact manifold $X$ whose derivatives up to order $k$
are in $L^{p}$. For any space $L^{p}_{k}$, we define its weight $w(k,p)=k-\frac{n}{p}$.

\textbf{Sobolev embedding theorem}
For two pairs $(p,k)$ and $(q,j)$, such that
$j\leq k$, $w(q,j)\leq w(p,k)$, there is a continuous inclusion:
$$L^{p}_{k}\rightarrow L^{q}_{j}$$

\textbf{Sobolev multiplication theorem}
If we have three pairs of numbers $(p,k)$, $(q,j)$ and $(r,i)$ such that
$i\leq \min\{k,j\}$, $w(r,i)\leq \min\{w(p,k),w(q,j)\}$ and
$w(r,i)\leq w(p,k)+w(q,j)$, we get the following bounded linear map:
$$L^{p}_{k}\times L^{q}_{j} \rightarrow L^{r}_{i}$$
given by pointwise product of functions.

\end{document}